\documentclass[12pt,amsfonts,amssymb,epsf,float]{article}
\usepackage{amsmath}
\usepackage{amssymb}
\usepackage{amsthm}
\usepackage{amsfonts}
\usepackage{graphicx}
\usepackage{enumerate}
\usepackage{bbm,bbold}
\usepackage{dsfont}

\usepackage{bbm,epsfig,graphics,epic,color,rotating,color}
\usepackage{ulem}

\numberwithin{equation}{section}

\include{(null)}

\def\re{{\mathrm Re\,}}
\def\im{{\mathrm Im\,}}

\newtheorem{theorem}{Theorem}[section]

\newtheorem{proposition}[theorem]{Proposition}

\newtheorem{definition}[theorem]{Definition}
\newtheorem{remark}[theorem]{Remark}

\theoremstyle{definition}

\usepackage{color}

\title{{Stochastic Dynamics of Einstein Matter-Radiation Model with Spikes}}

\author{Eugeny A. Pechersky, \\
Institute for Information Transmission Problems, \\
19, Bolshoj Karetny, Moscow, 127994, Russia \\
\textit{pech@iitp.ru}\\
\\
Anatoly A. Yambartsev \\
Instituto de Matem\'{a}tica e Estat\'{i}stica,
Universidade de S\~{a}o Paulo, \\
Rua do Mat\~{a}o, 1010, CEP 05508-090, S\~{a}o Paulo, SP, Brasil\\
\textit{yambar@ime.usp.br}\\
\\
Valentin A. Zagrebnov\\
D\'{e}partement de Math\'{e}matiques - AMU\\
Institut de Math\'{e}matiques de Marseille (UMR 7373)\\
CMI - Technop\^{o}le Ch\^{a}teau-Gombert - 39, rue F. Joliot-Curie \\
13453 Marseille Cedex 13, France\\
\textit{valentin.zagrebnov@univ-amu.fr}\\
}
\date{}

\begin{document}
\maketitle


\begin{abstract}
Motivated by the Einstein classical description of the matter-radiation dynamics we
revise a dynamical system producing spikes of the photon emission. Then we study the corresponding
stochastic model, which takes into account the randomness of spontaneous and stimulated atomic
transitions supporting by a pumping.
Our model reduces to Markovian density dependent processes for the inversion coefficient and
the photon specific value driven by small but fast jumps. We analyse the model in three limits:
the many-component, the mean-field, and the one-component. In the last case it is a positive
recurrent Markov chain with sound spikes.
\end{abstract}

\vspace{1.0cm}


\newpage
\tableofcontents


\section{Introduction}\label{sec.1}

We consider a stochastic model motivated by the Einstein description of the atom-radiation interaction
that takes into account the stimulated emission \cite{E17}. Our aim is to work out and to analyse a
two-component interacting Markov processes. The first component corresponds to the mean \textit{photon density}
$n$, whereas the second to the \textit{population inversion} coefficient $r$, which is equal to the mean
number of the ratio excited/disexcited two-level atoms.

In the nonstochastic limit these parameters are driving by a dynamical system with pumping that produces
excited atoms as well as the photon emission. We show that above a certain threshold for pumping and after
tuning of other parameters, this dynamical system manifests spikes of the photon density. In contrast to
attenuation in dynamical system, these spikes persists in our stochastic model. The aim of the paper is
to study this Markovian stochastic process with spikes in the three limit cases: a global many-atomic
Markov evolution, a one-atom evolution in random environment (mean-field approximation), and a free
one-unit stochastic evolution.

\section{Dynamical system with spikes of radiation}\label{sec.2}

\subsection{The Einstein matter-radiation equations }\label{sec.2.1}

We recall here main points of the Einstein matter-radiation theory with stimulated emission, see
\cite{E17} and \cite{Re} (Part II) for applications in laser theory.

To this aim let $N_1, N_2$ denote the numbers of two-level atoms, where $N_2$ corresponds to the number of \textit{excited} atoms and $N_1$ is the
number of \textit{non-excited} ones. (Recall that excited atoms are in the quantum state with higher
energy level $E_2 > E_1$ occupied. Then de-excitation of \textit{one} atom via transition:
$E_2 \rightarrow E_1$, creates \textit{one} photon.) By $r:= N_2/N_1$ we denote the
\textit{population inversion} coefficient. Let $\mathcal{N}\in\mathbb{Z}_+$ denote a \textit{total} number of photons in the system and
$n := \mathcal{N} / N_1 $ be their specific amount (\textit{density}).
Further, we denote by $p : = P/N_1 \geq 0$ the parameter of specific \textit{pumping} per atom
referring (as above for photons) to the number of non-excited atoms $N_1$.

The Einstein equations for the ``quantum" evolution system (QES), that we use for construction of our
stochastic model, are nonlinear ordinary differential equations for the number of excited atoms
$t \mapsto N_2 (t)$ and density of photons $t \mapsto n (t)$, for a \textit{fixed} $N_1$. These equations
have the form:
\begin{eqnarray*}
&&\Gamma \, \{\partial_t N_2\} = - A_{21}N_2 - w_{21}\, n N_2  + w_{12}\, n N_1 + P = \\
&&(spontaneous \ transitions : E_2 \rightarrow E_1 \ + \
stimulated \ transitions \ down:E_2 \rightarrow E_1 \ \\
&&+ \ stimulated \ transitions \ up:E_1 \rightarrow E_2 \ +
\ pumping \ to \ excite \ atoms) \\
&& \\
&&N_1 \, \{\partial_t n\} = {A_{21}}N_2 + w_{21}\, n N_2 - w_{12}\, n N_1 - b \, N_1 n =
(spontaneous \ emission \ +  \\
&&\ stimulated \ emission \ + \ lost \ by \ absorption \ exciting \ atom \ + \\
&&\ lost\ by \ leaking/radiation).
\end{eqnarray*}
Here we symbolically denote the rate of evolution by ``derivatives" $\{\partial_t N_2\}$ and
$\{\partial_t n\}$. Now there are few remarks in order \cite{Re}.

\textit{First}, all involved in these equations coefficients are non-negative. Moreover, the
Einstein transition coefficients are equal: $w_{21} = w_{12} > 0$, to ensure the \textit{detailed
balance} principle. Similarly to these coefficients, the spontaneous transition amplitude
$A_{21}> 0$ is entirely determined by the quantum properties of atoms. Note that the \textit{quality} of the
system is characterised by the ratio $\kappa = w_{21}/A_{21}$ and the leaking parameter $b$.
For quantum optical systems (laser, maser etc) one keeps $\kappa$ to be \textit{large}.

\textit{Second}, there are two \textit{external} parameters $b > 0$ and $p > 0$. The lost of photons due to
leaking/radiation depends e.g. on optical properties and geometry of the system varying $B$ in
a large scale. Without these precautions its value is of the \textit{order} of the Einstein transition
coefficients: $b \simeq w_{21} = w_{12}$.
The value of $p$ rules the rate of production of excited atoms. This is insured due to irradiation
of the system by external source of the pumping light with a higher frequency than the photons
of the system. This excites (in fact \textit{three-level}) atoms first to the level $E'_2 > E_2$
with a relatively large spontaneous transition amplitude $A'_{2 2}$ from $E'_2$ to the level $E_2$.
There they are living much longer because of the small $A_{21}$, waiting for enough value of
population inversion $r$ to produce an \textit{avalanche} of transitions down of de-exciting atoms
with a consequent \textit{spike} of the photon density in the system. The last is interpreted as
emission leaking out due to the term
$b N_1 n$.

\textit{Third}, it is clear that to realise this scenario the pumping $p$ should be strong enough to
make the population inversion coefficient $r$ to be of the order $\mathcal{O}(1)$.
\begin{remark}\label{rem:2.1}
The parameter $\Gamma > 1$ is usually \textit{large} \cite{Re}.
This reflects the fact that the \textit{time-scale} for evolution of $N_2 (t)$ for atoms
is (much) \textit{larger} than the \textit{time-scale} for evolution of $n(t)$ for photons.
For example, the frequency of oscillations of atoms are typically smaller then the light photon frequency.
Therefore, the atomic relaxation time is \textit{larger}), then the photon relaxation time.
\end{remark}
Taking into account these remarks we deduce equations of the matter-radiation Dynamical System (DS).
To this end we normalise the QES equations using as {parameter} $N_1$ and we choose the
Einstein transition coefficients as a reference: $w_{21} = w_{12} =1$.

Recall that the coefficient of leaking/radiation $b$ are usually of the {same} order
of magnitude as transition coefficients. It is convenient to put $b=:\beta^{-1}> 0$. We also define
$\alpha = A_{21}$, for (small) values of non-negative parameter $\alpha \geq 0$. Note that the
system with the infinitely hight \textit{quality} corresponds to the limit case $\alpha = 0$, when
the \textit{spontaneous} transitions/emission are completely suppressed.

Then DS equations take the form
\begin{eqnarray}\label{DS1}
\Gamma\partial_t r &=& (- \alpha r -  n r + n) + p  \ ,  \\
\partial_t n &=& (\alpha r  +  n r -  n) - \beta^{-1}n \ . \label{DS2}
\end{eqnarray}

\subsection{Dynamical system: stationary points and linearisation}\label{sec.2.2}

Consider the stationary solutions of the system (\ref{DS1}), (\ref{DS2}):
\begin{equation}
\left\{ \begin{array}{rcl}
0 &=& (- \alpha r - n r + n) + p \ ; \\
0 &=& (\alpha r  + n r -  n) - \beta^{-1}n \ .
\end{array} \right.
\end{equation}
It provides only one stationary solution:
\begin{equation}\label{s2}
\left\{ \begin{array}{rcl}
r^* & = & p(1 + \beta)/(\beta p + \alpha) \ , \\
n^* &=& \beta p \ .
\end{array} \right.
\end{equation}

Note that for large $p$, or small $\alpha$, one gets $r^* = \mathcal{O}(1)$, which confirms the
relevance of the Dynamical System (\ref{DS1}), (\ref{DS2})) in the regime of the high population
inversion coefficient, when the spikes are producing.

Linearisation: $r = r^* + \delta r, n= n^* + \delta n $, of the
system (\ref{DS1}), (\ref{DS2}) in the vicinity of the stationary point $(r^*, n^*)$
gives the time evolution equation:
$$
\partial_t \binom{\delta r}{\delta n} = \left[ \begin{array}{cc} - (\alpha + n^*)/{\Gamma} &
(-r^* + 1)/{\Gamma} \\  \alpha + n^* & r^* - (1+ \beta^{-1})\end{array} \right]
\binom{\delta r}{\delta n}
\ .
$$
Then by (\ref{s2}) the corresponding characteristic equation gets the form
$$
0 = \det \left[ \begin{array}{cc} - (\beta p + \alpha)/{\Gamma} -\lambda & -(p -\alpha)/{(\beta p +\alpha)\Gamma} \\
\beta p + \alpha & -\alpha(1+\beta)/(\beta p +\alpha)\beta -\lambda \end{array} \right]
=
$$
$$
= \lambda^2 + \lambda \, \left[\frac{1}{\Gamma}(\beta p + \alpha) + \frac{\alpha(1+\beta)}
{(\beta p +\alpha)\beta}\right] + \frac{1}{\Gamma}\frac{\beta p +\alpha}{\beta} \ .
$$
Let us introduce $z: = \beta p + \alpha$. Then two roots of the characteristic equation are
\begin{eqnarray}\label{roots}
\lambda_{1,2} &=& - \frac{1}{2}\left[\frac{z}{\Gamma}\, + \frac{\alpha(1+\beta)}{z\beta}\right]
\pm \sqrt{\frac{1}{4}\left[\frac{z}{\Gamma} + \frac{\alpha(1+\beta)}{z\beta}\right]^2 - \frac{z}{\Gamma\beta}}
\\ &=& - \frac{1}{2}\left[\frac{z}{\Gamma}\, + \frac{\alpha(1+\beta)}{z\beta}\right] \pm \sqrt{\Delta}.
\nonumber
\end{eqnarray}

Now we observe that if $\Delta \ge 0$,  then $\lambda_{1,2} < 0$. This  means that the fixed point
$(r^*, n^*)$ is a \textit{stable node}. On the other hand, for $\Delta < 0$ the fixed point
$(r^*, n^*)$ is a \textit{stable focus}, when $\im{\lambda_{1,2}} \neq 0$ and $\re{\lambda_{1,2}} <0$.
The boundary between two regimes is defined by equation $\Delta = 0$, or
\begin{equation}\label{discrim}
\frac{1}{4}\left[\frac{z}{\Gamma} + \frac{\alpha(1+\beta)}{z\beta}\right]^2 - \frac{z}{\Gamma\beta} = 0 \
\end{equation}

First we consider the case of the system with the \textit{infinitely hight} quality: $\alpha = 0$.
Then for fixed $p, \beta > 0$ the boundary value for $\Gamma$ is $\Gamma_{0} = \beta^{2} p/4$. By
(\ref{discrim}) one gets $\Delta (\Gamma \leq \Gamma_{0}) \geq  0$. For these values of $\Gamma $
the point $(r^*, n^*)$ is a {stable node}. As far as  concerns two roots (\ref{roots}) of the
corresponding to $\alpha = 0$ characteristic equation, one gets
$\lim_{\Gamma \rightarrow 0} \lambda_{1}^{0}(\Gamma) = - \infty$,
whereas $\lim_{\Gamma \rightarrow 0} \lambda_{2}^{0}(\Gamma) = - 1/\beta$ and
$\lim_{\Gamma \rightarrow \Gamma_{0}} \lambda_{1,2}^{0}(\Gamma) = - 2/\beta$. Similarly, (\ref{discrim})
yields that $\Delta (\Gamma > \Gamma_{0}) < 0$. For these values of $\Gamma $ the point
$(r^*, n^*)$ is a {stable focus}. Note that in this case
$\lim_{\Gamma \rightarrow \infty} \lambda_{1,2}^{0}(\Gamma) = 0$.
We conclude by remark that the frequency of rotations (oscillations):
$\im{\lambda_{1}^{0}(\Gamma)} = \sqrt{- \Delta (\Gamma)}$, in the stable focus monotonously increases
for $\Gamma > \Gamma_{0}$ and reaches the maximum $\im{\lambda_{1}^{0}(\Gamma_{0}^{*})} = 1/\beta$ at
$\Gamma_{0}^{*} = \beta^{2} p/2$. Then $\im{\lambda_{1}^{0}(\Gamma)}$ monotonously decreases to zero
for $\Gamma \rightarrow \infty$.

Now we consider the case $\alpha > 0$. Then by (\ref{discrim}) for fixed $\alpha, p, \beta > 0$ the
boundary values for $\Gamma$ are
\begin{equation}\label{gamma}
\Gamma_{1,2} = \frac{\beta^2z^2}{\alpha^2(1+\beta)^2} \left[ \frac{2z -
\alpha(1+\beta) }{\beta} \mp 2\sqrt{ \frac{z}{\beta} (p - \alpha) } \right] \ .
\end{equation}

Note that for $\alpha \geq p$ the complex $\Gamma_{1,2}$ that discriminant $\Delta(\Gamma) > 0$
for all $\Gamma \geq 0$. Therefore, by (\ref{roots}) the point $(r^*, n^*)$ is always a {stable node}.

For $\alpha < p$ the two real $\Gamma_{1,2}$ reflect the following behaviour of $\Delta(\Gamma)$:\\
(a) On the interval $[0, \Gamma^*]$ it is monotonously decreasing from $+\infty$ and reaches the
minimum $\Delta(\Gamma^*) = - (\beta p - \alpha)/z\beta^{2} < 0$ at
$\Gamma^* = \beta z^2 /(2 \beta p + \alpha (1-\beta))$.\\
(b) On the interval $[\Gamma^*, \infty)$ discriminant $\Delta(\Gamma)$ is monotonously increasing
to the limit value $[{\alpha(1+\beta)}/2{z\beta}]^2 > 0$.

Therefore, by (\ref{roots}) on gets that on the interval $[0, \Gamma^*]$  the point $(r^*, n^*)$  changes at
$\Gamma_{1}$ from the {stable node} to a stable focus with the maximal frequency
$\im{\lambda_{1}(\Gamma^*)} = \sqrt{- \Delta (\Gamma^*)}$. Then on the interval $[\Gamma^*, \infty)$
the point $(r^*, n^*)$  changes at $\Gamma_{2}$ from the {stable focus} to a stable node.
Note that the \textit{high} pumping limit: $p/\alpha \rightarrow \infty$, yields
$\lim_{p/\alpha \rightarrow \infty} \Gamma_{1} = \Gamma_0$ and
$\lim_{p/\alpha \rightarrow \infty} \Gamma_{2} = +\infty$. This corresponds to the case of only one
transition point at $\Gamma_0$, that we considered above.

Illustration of behaviour of DS (\ref{DS1}), (\ref{DS2}) is presented in Figures \ref{fig1} and
\ref{fig2}, where $r = \rho$ and $n = \nu$. In agreement with our analysis one observes a sound
\textit{spike} of the photon density $\nu$ for $\Gamma = 100$, but not for $\Gamma = 2$.

\section{Stochastic Fluctuations of Dynamical System}\label{sec.3}

The inference of dynamical system equations (\ref{DS1}), (\ref{DS2}) was based on
the Einstein model of the matter-radiation interaction. They yield a time evolution for the
mean densities (specific numbers) of photons and the population inversion coefficient for
excited/non-excited atoms. In the framework of this, in fact \textit{classical} model, there are
two ways to include \textit{quantum} evolution into dynamical system (\ref{DS1}),(\ref{DS2}).

One way is to return back to the basic quantum mechanical equations.
Another way is to mimic this evolution by considering probabilistic jumps between excited/non-excited
atomic states coupled to stochastic emissions of photons as a result of two Markov processes.
In this way the "quantum" evolution can be treated in the spirit of Einstein's classical description
of the matter-radiation interaction.

Below we provide a construction of such coupled Markov processes for the random population inversion
coefficient and the photon density.

\subsection{Markov process driven by dynamical system}\label{sec.3.1}

Note that (DS) Dynamical System (\ref{DS1}), (\ref{DS2}) describes evolution of \textit{collective}
(intensive) variables: the coefficient of inversion $r = N_2 / N_1$ and the density of photons
$n = \mathcal{N} / N_1$.

To establish a link between DS and the corresponding Markov process we first take into account
a difference of the time-scale for evolution of $t \mapsto r (t)$ and $t \mapsto n (t)$, where
$t\geq 0$. Then for a time-increment $\Delta t$ and the corresponding $t^* \in [t, t+ \Delta t]$ we
obtain from (\ref{DS1}), (\ref{DS2}) the representation:
\begin{eqnarray}\label{DS12}
r(t+ \Delta t)  &=& r(t) + \frac{\Delta t}{\Gamma}
[(- \alpha r -  n r + n)(t^*) + p]  \ ,  \\
n(t+ \Delta t)  &=& n(t) + \Delta t
[(\alpha r + n r - n)(t^*) - \beta^{-1}n(t^*)] \ . \label{DS22}
\end{eqnarray}
According to (\ref{DS12}) and (\ref{DS22}) the time-scale difference imposed
by $\Gamma$ in (\ref{DS1}, \ref{DS2}) can be
transformed into scale difference of \textit{increments} for a \textit{random} version of
functions $r, n$.
\begin{remark}\label{rem:3.1}
To proceed further we fix the \textit{normalising} size parameter and identify it with
the time-independent number of non-excited atoms $N: = N_1 \in \mathds{N}$, see Section \ref{sec.2.1}.
Therefore, the total number of atoms is \textit{not} fixed, but it plays no role in evolution.
Then we rename variables, which are \textit{specific values} normalised to number $N$ of the non-excited
atoms. We denote them by $r_N := N_2 / N$ (the coefficient of inversion, or density of excited atoms) and
$n_N := \mathcal{N} / N$ (density of photons). Then we keep notations: $r = \lim_{N\rightarrow\infty}r_N$
and $n = \lim_{N\rightarrow\infty}n_N$ for their limits.
\end{remark}


Now we note that according to Einstein equations (Section \ref{sec.2.1}) the \textit{elementary}
transition processes in our system are: \\
- the \textit{one-atom} excitation/de-excitation accompanied by \textit{one-photon} absorption/creation, \\
- the \textit{one-atom} excitation due to pumping $p = P/N$, and \\
- the \textit{one-photon} leaking/radiation with the rate $b$. \\
Then by the size normalisation and by taking into account
different time-scale (\ref{DS12}), (\ref{DS22}) the variables $r_N, n_N$ jump
(after the \textit{sojourn time}) with increments ${1}/{N\Gamma}, {1}/{N}$ either \\
as $(r_N, n_N)\rightarrow (r_N \pm 1/N\Gamma , n_N \mp 1/N)$, or \\
as $(r_N, n_N)\rightarrow (r_N + 1/N\Gamma , n_N)$ and
$(r_N, n_N)\rightarrow (r_N , n_N - 1/N)$. \\
This means that variable $r_N$ is jumping on the lattice $\frac{1}{\Gamma N}\mathbb{Z}_+ \subset \mathbb{R}_+$,
whereas variable $n_N$ jumps on the lattice $\frac{1}{N}\mathbb{Z}_+ \subset \mathbb{R}_+$.

Summarising the elementary transitions and the choice of the corresponding rates defined by the
Einstein QES, or DS (\ref{DS1}), (\ref{DS2}) via densities $r_N$, $n_N$, we obtain the following list of
Markovian jumps and intensities:
\begin{enumerate}
\item[(1)] with the rate $N r_N n_N$ it occurs $(r_N,n_N) \to (r_N-\frac{1}{\Gamma N}, n_N + \frac{1}{N})$, i.e.
the number of exited atoms decreases by one and creates one photon; the rate depends on the product $r_N n_N$,
which is interpreted
as transition with stimulated emission since it is proportional to the density of photons $n_N$ in the system;

\item[(2)] with the rate $\alpha N r_N$ it occurs $(r_N,n_N) \to (r_N-\frac{1}{\Gamma N}, n_N+\frac{1}{N})$, i.e.
the number of exited atoms decreases by one with the rate proportional to the density of excited atoms $r_N$
(emission due to spontaneous transition);

\item[(3)] with the rate $N n_N$ it occurs $(r_N,n_N) \to (r_N+\frac{1}{\Gamma N}, n_N-\frac{1}{N})$, i.e.
excitation of atom with absorption of one photon (stimulated transition up).

\item[(4)] with the rate $N p$ it occurs $(r_N,n_N)  \to (r_N+\frac{1}{\Gamma N}, n_N)$, where $p\in \mathbb{R}_+$
 is the pumping parameter defined before; it means that an atom was exited without changing the number of photons
 in the system (pumping of excited atoms);

\item[(5)] with the rate $\beta^{-1}N n_N$ it occurs $(r_N,n_N) \to (r_N, n_N - 1/N)$, i.e. lost of photons via
leaking by radiation.
\end{enumerate}
\begin{remark}\label{rem:3.2}
Note that our choices in (1)- (5) corresponds to the stochastic evolution of a large many-component
($N$-units) system which evolves via many \textit{small} individual jumps of order $1/N\Gamma$ and $1/N$,
but at a \textit{fast} rate of order $N$. In fact, this is an appropriate mathematical background corresponding
to the physical nature of evolution of the Einstein matter-radiation model, which
in nonrandom ``smooth" approximation is DS (\ref{DS1}), (\ref{DS2}).

Initiated by Kurtz in seventies \cite{K}, \cite{EK}, this class of Markov processes in the
limit of large $N$ is known under the \textit{law large number scaling} \cite{EK} as well as  the
\textit{mean-field approximation} \cite{Kol} or the \textit{fluid limit} \cite{DN}. The motivation
is that for large $N$ the density dependent Markov process with small increments of $O(1/N)$ but
with intensities of order $N$ may be approximated by a nonrandom trajectory verifying a differential
equation constructed the rate intensities \cite{K}.

We also note that if symmetric jumps have the size of order $1/\sqrt{N}$  still with the rate
of order $N$, then this approximation for the Markov processes is called \textit{central limit}, or
\textit{diffusive} scaling, see \cite{EK} and \cite{JS}.
\end{remark}

In the present paper we move (in a certain sense) \textit{backward}: starting from the Einstein QES and DS
(\ref{DS1}), (\ref{DS2}) for \textit{densities} we construct a driven by DS  Markov jump process with
small increments of the order $O(1/N)$ motivated by the fast quantum transitions in the $N$-component system,
where $N$ is of the order of number of atoms.

\smallskip

To proceed we recall first some key notations, definition and statement due to \cite{EK}, \cite{K}.
They are indispensable for our construction and analysis of the random process generated by the Einstein
QES and DS (\ref{DS1}), (\ref{DS2}) for physical \textit{specific values}, i.e. \textit{densities} of
the excited atoms and the photon number $(r_N, n_N)$.
\begin{definition}\label{def:3.1}
A family of continuous-time jump Markov process $\{\Phi_{N}(t)\}_{t\geq 0}$, parameterised by
$N \in \mathbb{N}$,  with the space of states $S \subseteq \mathbb{Z}^{d}$ is called
\textit{density dependent} if the corresponding transitions: $s_1 \rightarrow s_2$, for any
states $s_1, s_2 \in S$, have the rates $T_{N}(s_1, s_2) \geq 0 $, which have the
following form
\begin{equation}\label{eq:3.1}
T_{N}(s_1, s_2) := N \beta_{s_2 - s_1} (s_1 /N)  \ , \ s_1, s_2 \in S \ .
\end{equation}
The set of functions $\{E \ni x \mapsto \beta_{l}(x)\}_{l \in S}$ are defined and non-negative
on a subset $E \subset \mathbb{R}^{d}$.
\end{definition}

Let $\{\nu_{l\in S}(t)\}_{t\geq 0}$ be family of \textit{independent} \textit{unit-rate} Poisson processes
that count the occurrences of the events corresponding to jumps of the process $\Phi_{N}(t)$ by $l\in S$.
Then by virtue of (\ref{eq:3.1}) the jump Markov process $\{\Phi_{N}(t)\}_{t\geq 0}$ verifies the
stochastic equation:
\begin{equation}\label{eq:3.2}
\Phi_{N}(t) = \Phi_{N}(0) +
\sum_{l\in S} \, l \, \nu_{l}\left(\int_{0}^{t} d\tau \, N \beta_{l} (\Phi_{N}(\tau)/N)\right) \ .
\end{equation}

To study this class of Markov chains we follow \cite{EK}, \cite{K} and rescale $\Phi_{N}(t)$ into the
\textit{density} process $\{\phi_{N}(t) = \Phi_{N}(t)/N\}_{t\geq 0}$, $N \in \mathbb{N}$.
Then (\ref{eq:3.2}) yields the stochastic equation for density processes
\begin{equation}\label{eq:3.3}
\phi_{N}(t) = \phi_{N}(0) +
\frac{1}{N}\sum_{l\in S} \, l \, \int_{0}^{t} d\tau \, N \beta_{l} (\phi_{N}(\tau))
+ \frac{1}{N}\sum_{l\in S} \, l \,
\widehat{\nu}_{l}\left(\int_{0}^{t} d\tau \, N \beta_{l} (\phi_{N}(\tau))\right).
\end{equation}
Here $\{\widehat{\nu}_{l}(t) : = {\nu}_{l}(t) - t\}_{t\geq 0}$ are \textit{compensated} unit-rate
Poisson processes. Note that $\phi_{N}\in (1/N)\mathbb{Z}^{d}$, i.e., the jumps of the process
(\ref{eq:3.3}) have increments $O(1/N)$ whereas the intensities (\ref{eq:3.1}) are of the order $O(N)$.
The corresponding to the process (\ref{eq:3.3}) generator has the form
\begin{equation}\label{eq:3.4}
(\mathcal{L}_{N}\mathcal{G})(x) =  N \sum_{l\in S}\beta_{l}(x)\ [\mathcal{G}(x + l/N) - \mathcal{G}(x)] \ , \
x \in (1/N)\mathbb{Z}^{d} \ ,
\end{equation}
for functions $\mathcal{G}$ with compact supports on the lattice $(1/N)\mathbb{Z}^{d}$.

We note that these observations suggest that for increasing parameter $N$ the c\`{a}dl\`{a}g trajectories
of the process (\ref{eq:3.3}) approximate a continuous nonrandom trajectory $\phi(t)$.
Indeed, let us take into account that the Law of Large Numbers (LLN) for the compensated Poisson process
$\{\widehat{\nu}_{l}(t)\}_{t\geq 0}$ implies for each $t_0 > 0$  and $f(t) \geq 0$
\begin{equation}\label{eq:3.5}
a.s.-\lim_{N \rightarrow \infty} \sup_{t \leq t_0} \ \frac{1}{N} \ \widehat{\nu}_{l}(N f(t)) = 0 \ ,
\end{equation}
in the representation (\ref{eq:3.3}). Then the density process $\{\phi_{N}(t)\}_{t\geq 0}$ converges
to deterministic solution $\phi(t)$ of differential equation
\begin{equation}\label{eq:3.6}
\partial_t \phi(t) = \sum_{l\in S} \, l \, \beta_{l} (\phi(t)) \ ,
\end{equation}
which express the LLN for this kind of process. Now we formulate the exact statement, which is due to
\cite{K}, \cite{EK}:
\begin{proposition}\label{prop:3.1}
Let $\sum_{l\in S} \, |l| \, \sup_{\phi \in K}\beta_{l} (\phi) < \infty$ for each compact
$K \subseteq E \subset \mathbb{R}^{d}$, and the function $\phi \mapsto \sum_{l\in S} \, l \, \beta_{l} (\phi)$
be Lipschitz continuous. If the Markov process $\{\phi_{N}(t)\}_{t\geq 0}$ satisfies the stochastic equation
(\ref{eq:3.3}) with initial condition such that $a.s.-lim_{N\rightarrow \infty}\phi_{N}(0) = \phi(0)$ and
$\phi(t)$ is the corresponding solution of the Cauchy problem (\ref{eq:3.6}), then for any $t \geq 0$
one gets
\begin{equation}\label{eq:3.7}
a.s.-\lim_{N \rightarrow \infty} \sup_{\tau \leq t} |\phi_{N}(\tau) - \phi(\tau)| = 0 \ .
\end{equation}
\end{proposition}
We conclude by remark concerning the time-scale change: $\partial_t  \rightarrow \Gamma \partial_t $.
Then (\ref{eq:3.6}) transforms into equation
\begin{equation}
\partial_t \phi_{\Gamma}(t) = \sum_{l\in S} \, \frac{l}{\Gamma} \ \beta_{l} (\phi_{\Gamma}(t)) =
\sum_{l/\Gamma \in S_{\Gamma}} \, \frac{l}{\Gamma} \ \widehat{\beta}_{l/\Gamma} (\phi_{\Gamma}(t))\ ,
\end{equation}
where $S_{\Gamma} \subseteq (1/\Gamma)\mathbb{Z}^{d}$. Therefore, the time-scale change is equivalent
to the corresponding change of increment of jumps that we found in (\ref{DS12}) and (\ref{DS22}).

\subsection{{Global Markov process driven by random densities}}\label{sec.3.2}
Following our \textit{backward} strategy we first establish  the analogue of (\ref{eq:3.6}) for the
Einstein QES and equations (\ref{DS1}), (\ref{DS2}). To this end we rewrite them as the system of
differential equations for \textit{two}-component density $X(t)$:
\begin{equation}\label{s1}
\partial_t X(t) = \partial_t \binom{{r}}{{n}}(t) = \sum_{l\in S} \, l \, \beta_{l} (X(t)) \ ,
\end{equation}
with initial \textit{non-negative} conditions: $r(0)=r_0$, $n(0)=n_0$. Here the \textit{basic} set
of the \textit{vector-valued} increments: $\{l_j\}_{j=1}^4 \subset S \subseteq ({1}/{\Gamma})\mathbb{Z}_+
\oplus \mathbb{Z}_+ $ , and the corresponding intensities $\{\beta_{l_j}\}_{j=1}^4$ are such that
\begin{equation}\label{eq:3.8}
\sum_{\{l_j : j = 1, \ldots, 4\}} l_j \, \beta_{l_j}  = \binom{-1/\Gamma}{1} \ (r n + \alpha r )
 + \binom{1/\Gamma}{-1} \ n + \binom{1/\Gamma}{0} \  p + \binom{0}{-1} \ \beta^{-1} n \ .
\end{equation}
Then taking into account representation (\ref{eq:3.3}) and (\ref{eq:3.8}) one gets the stochastic
equation for two-component Markov process $\{X_{N}(t)\}_{t\geq 0}$ driving random densities for a
given parameter $N \in \mathbb{N}$:
\begin{equation}\label{eq:3.9}
X_{N}(t) = X_{N}(0) +
\sum_{\{l_j : j = 1, \ldots, 4\}} \, \frac{l_j}{N} \  \nu_{l_j}\left(\int_{0}^{t} d\tau \,
N \beta_{l_j} (X_{N}(\tau))\right) \, , \ X_{N}(t) := \binom{{r_{N}(t)}}{{n_{N}(t)}} \ .
\end{equation}
Note that the basic increments of the process (\ref{eq:3.9}) are: $l_j/N \in ({1}/{\Gamma}N)\mathbb{Z}_+
\oplus ({1}/N)\mathbb{Z}_+ $, for $j= 1,2,3,4$.
By construction and by explicit values of $\{\beta_{l_j}\}_{j=1}^4$, the process (\ref{eq:3.9}) verifies
conditions of Proposition \ref{prop:3.1}. Therefore, we obtain for (\ref{s1}) and (\ref{eq:3.9}) the LLN
in the form:
\begin{equation}\label{eq:3.10}
a.s.-\lim_{N \rightarrow \infty} \sup_{\tau \leq t} |X_{N}(\tau) - X(\tau)| = 0 \ ,
\end{equation}
for any $t \geq 0$.
\begin{figure}[ht!]
\begin{center}
\includegraphics[width=0.8\linewidth]{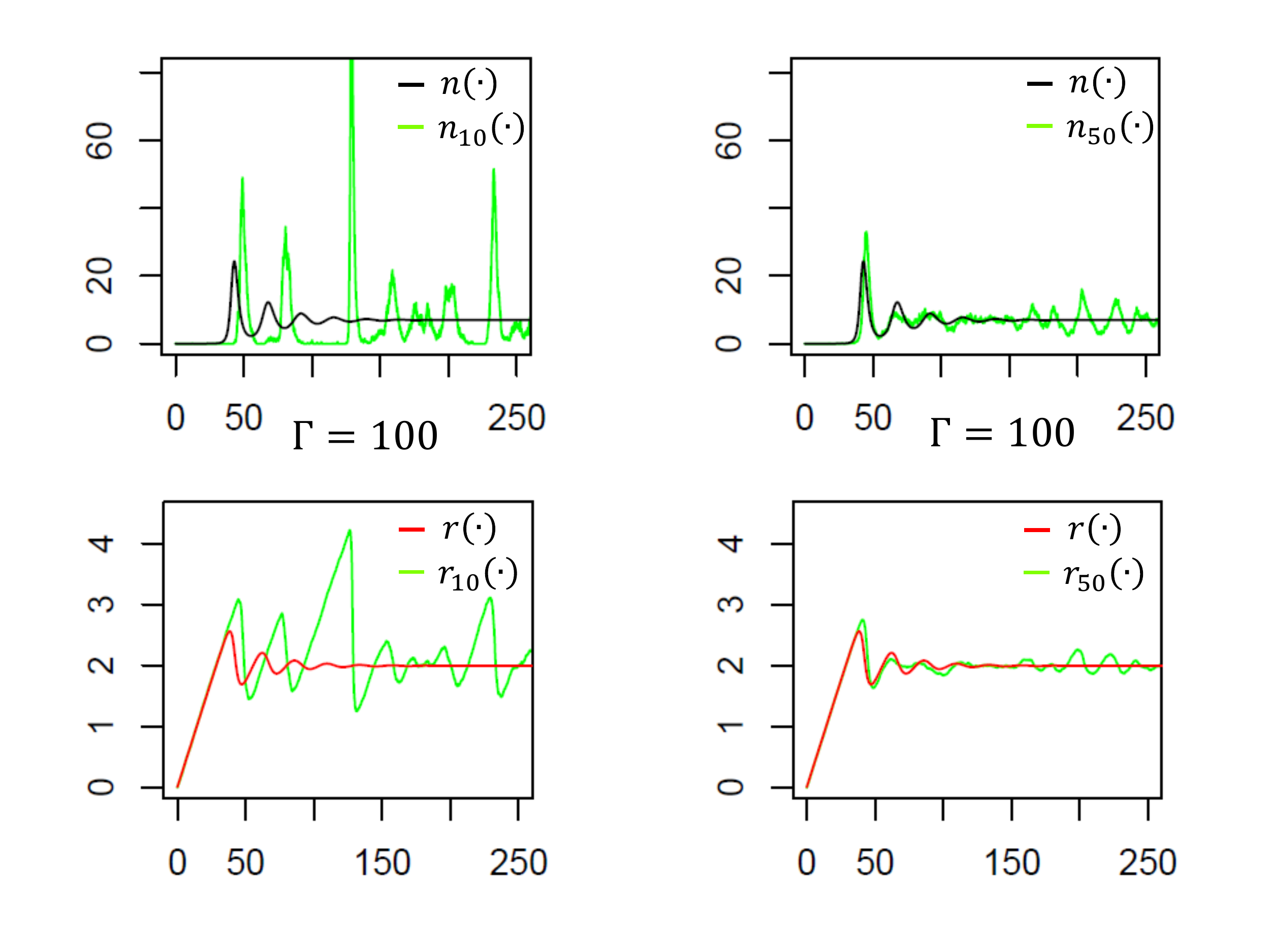}
\end{center}
\caption {Dynamical system for $\Gamma=100$ : on the top the solution for $n$ (black), on the
bottom for $r$ (red). Both starting at $0,01$.
Trajectory samples of stochastic global Markov processes (green): for $N=10$ (top-left, bottom-left)
and for $N=50$ (top-right, bottom-right), see definition (\ref{eq:3.9}). The chosen parameters are:
$\alpha=0,01, \Gamma=100, P=7$ and $\beta=1$.}%
\label{fig1}
\end{figure}
Now taking into account the list (1)-(5), Section \ref{sec.3.1}, for independent elementary transitions
on $S$ with the corresponding intensities for the two-component jump Markov process (\ref{eq:3.8}), we can
reconstruct the right-hand side of (\ref{eq:3.4}) for generator $\mathbf{L}_N$ of \textit{correlated}
two-component Markov process:
\begin{eqnarray}\label{Mark-gen}
(\mathbf{L}_N G)\left(x,y\right) &=& N x y  [G(x-\frac{1}{\Gamma N}, y + \frac{1}{N})-G(x, y)] \\
&+& N \alpha x [G(x -\frac{1}{\Gamma N}, y+\frac{1}{ N})-G(x,y)] \nonumber \\
&+& N y [G(x+\frac{1}{\Gamma N}, y-\frac{1}{ N})-G(x,y)] \nonumber \\
&+& N p [G(x+\frac{1}{\Gamma N}, y)-G(x,y)] \nonumber \\
&+& N \beta^{-1} y [G(x, y - \frac{1}{N})-G(x,y)] \ . \nonumber
\end{eqnarray}
Here bounded functions $G$ have compact supports on the lattice $\Lambda_N$, where
$(x,y) \in \Lambda_N :=({1}/{\Gamma N})\mathbb{Z}_+ \times ({1}/{N})\mathbb{Z}_+ $.

\begin{figure}[ht!]
\begin{center}
\includegraphics[width=0.8\linewidth]{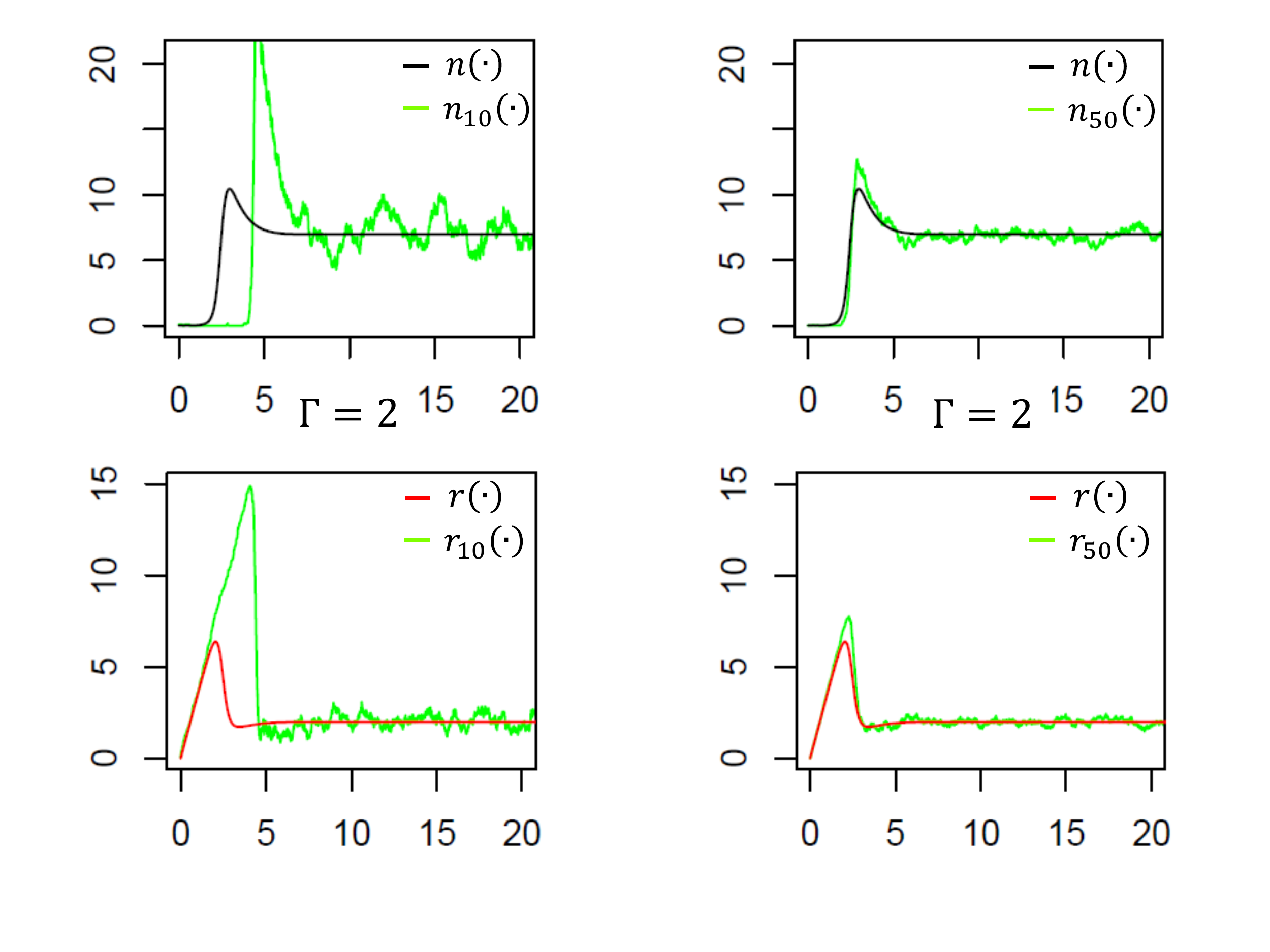}
\end{center}
\caption {Dynamical system for $\Gamma=2$ : on the top the solution for $n$ (black), on the bottom
for $r$ (red).
Both starting at $0,01$.
Trajectory samples of stochastic global Markov processes (green): for $N=10$ (top-left, bottom-left)
and for $N=50$ (top-right, bottom-right), see definition (\ref{eq:3.9}).
The chosen parameters are: $\alpha=0,01, \Gamma=2, P=7$ and $\beta=1$. }%
\label{fig2}
\end{figure}

Behaviour of the \textit{empirical} random density $X_N(t)$ for finite number of units $N$ is
illustrated by Figures \ref{fig1} and \ref{fig2}.
\begin{remark}\label{rem:3.2.1}
Note that coefficients of transition intensities in $\mathbf{L}_N$ (\ref{Mark-gen}) are taken
\textit{not} from the driving dynamical system (\ref{s1}) for densities, but they are generating
recursively, i.e., step-by-step along trajectories corresponding to the algorithm given by (\ref{Mark-gen}).

We also note that for jumps on $\Lambda_N =({1}/{\Gamma N})\mathbb{Z}_+ \times ({1}/{N})\mathbb{Z}_+ $
the values of components of the Markov process $X_N(t)$ (\ref{eq:3.9}) are unbounded for any $N \geq 1$.
\end{remark}
As we mentioned above the process $X_N(t)$ converges to a differentiable trajectory
of DS (\ref{DS1}), (\ref{DS2}).
The illustration of this statement is visible in Figures \ref{fig1} ($\Gamma=100$) and Figure \ref{fig2} ($\Gamma=2$),
when one compares the realisations of $X_N(t)$ for $N=10$ and for $N=50$. The green trajectories get
closer to the limit (\ref{eq:3.10}) for increasing $N$. Note that the spikes of the photon number for
$n_{N=10}$ are more sound than in the case of $n_{N=50}$. Similarly, the fluctuations of the coefficient of
inversion $r_{N=10}$ are more visible than those of $r_{N=50}$. Recall that these fluctuations are caused by
fluctuations of the excited atoms with respect to the fixed number of non-excited atoms $N$.

Using this interpretation one can represent the total random number of excited atoms: $\xi_N(t) = N r_N(t)$,
and the random number of photons: $\eta_N(t) = N n_N(t)$, as
\begin{equation}\label{eq:3.11}
\xi_N(t) := \sum_{j=1}^N \xi(t,j) \ \ {\rm{and}} \ \ \eta_N(t) := \sum_{j=1}^N \eta(t,j) \ .
\end{equation}
Here for each index $j$: $1 \leq j \leq N$, we introduce trajectories of identical, independent, compound
\textit{unit-rate} Poisson (\textit{telegraph}) processes: $\{\xi(t,j)\}_{t\geq 0}$  with values
$\xi(t,j) \in \{0,1\}$ for the jump-increment $1/\Gamma$, and $\{\eta(t,j)\}_{t\geq 0}$ with values $\eta(t,j)
\in \{0,1\}$ for jump-increments equal to $1$.

Therefore, by (\ref{eq:3.11}) the density Markov process $X_N(t)$ (\ref{eq:3.9}) corresponds to empirical
arithmetic means over \textit{trajectories}: $\{\xi(t,j)\}_{t\geq 0, j=1,\ldots,N}$ and
$\{\eta(t,j)\}_{t\geq 0, j=1,\ldots,N}$. Then the limit (\ref{eq:3.10}) expresses the LLN for the
arithmetic means of (\ref{eq:3.11}) over parameter $N$ of the number of non-excited atoms.

\subsection{Single Markov trajectories in a mean-field approximation}\label{sec.3.3}
In this section instead of arithmetic means (\ref{eq:3.11}) corresponding to \textit{collective} random
variables $r_N(t)$ and $n_N(t)$ we consider the \textit{individual} trajectories
$\{\xi(t,j)\}_{t\geq 0}$ and $\{\eta(t,j)\}_{t\geq 0}$.

First we note that by virtue of (\ref{Mark-gen}) and the algorithm of calculations of transition
intensities in Remark \ref{rem:3.2.1} besides the evident correlation between $\xi(t,j)$ and $\eta(t,j)$
for the same $j$, this trajectories are correlated since (\ref{Mark-gen}) is not a sum of independent
generators $\{\mathbf{L}_{N}^{(j)}\}_{j=1,\ldots,N}$.
To take into account the impact of these correlations on a single Markov trajectory $j$:
\begin{equation}\label{eq:3.12}
 X_{N}^{(j)}(t) := \binom{{\xi(t,j)}}{{\eta(t,j)}} \ ,
\end{equation}
we use a  \textit{mean-field approximation}. This approximation splits the above mentioned correlations
and allows to present generator (\ref{Mark-gen}) as a sum generators for single trajectories correlated
only in the mean-field.

To this aim we recall that coefficients of transitions in generator (\ref{Mark-gen}) are related to the
instant values of density processes $x= r_N$ and $y = n_N$. Then one can use the representations:
$N r_N n_N = \frac{1}{2} \sum_{j=1}^{N} ( x_{j} n_N + r_N y_{j})$,
$N r_N :=  \sum_{j =1}^{N}  x_{j}$, $N n_N :=  \sum_{j =1}^{N}  y_{j}$,
to rewrite $\mathbf{L}_N$ (\ref{Mark-gen}) identically in the form:
\begin{equation*}
\mathbf{L}_N  := \sum_{j =1}^{N} \  \mathbf{\mathcal{L}}_{j} \ ,
\end{equation*}
where each $\mathbf{\mathcal{L}}_{j}$ represents generator of a \textit{singe} trajectory correlated to
others. In the \textit{mean-field} approximation (which includes $n_N \approx n$ and $r_N \approx r$)
we define generator $\mathbf{\mathcal{L}_{m-f}}$ for the single-trajectory (\ref{eq:3.12}) as:
\begin{eqnarray}\label{Mark-gen-single}
(\mathbf{\mathcal{L}_{m-f}} \ g)\left( x,y \right) &:=& \frac{1}{2}(x \, n +
r \, y)[g(x -\frac{1}{\Gamma}, y + 1)- g(x,y)] \\
&+& \alpha x [g(x -\frac{1}{\Gamma}, y + 1)- g(x,y)] \nonumber \\
&+& y [g(x+\frac{1}{\Gamma}, y - 1)- g(x,y)] \nonumber \\
&+& p [g(x +\frac{1}{\Gamma}, y)- g(x,y)] \nonumber \\
&+& \beta^{-1} y [g(x, y - 1)- g(x,y)] \ , \nonumber
\end{eqnarray}
for functions $g$ with compact supports on the lattice $\Lambda:=({1}/{\Gamma})\mathbb{Z}_+
\times \mathbb{Z}_+ $. Note that by equivalence of trajectories we skipped in (\ref{Mark-gen-single})
the index $j$, and that parameters $(r,n)$ are driven by the corresponding dynamical system for
densities (\ref{s1}). {Let denote a Markov process governed by the generator (\ref{Mark-gen-single})
by $Z(t)=(\xi(t),\eta(t))$.}

Functions $r=r(t)$ and $n=n(t)$ appearing in the expression for the first rate in
(\ref{Mark-gen-single}) depend on time following the evolution
(\ref{s1}). But it is known that this dynamical system converges rapidly to the stationary point
$(r^* , n^*)$ (\ref{s2}). Thus, in order to analyse the behaviour
of the corresponding process it is reasonable to consider the Markov process, when the transition rates
depend only on $(r^* , n^*)$. It makes the process homogeneous in time and facilitate the analysis.
Denote such process by $Z_0(t)=(\xi_0(t),\eta_0(t))$. The following theorem states the ergodicity of
the process.
\begin{theorem}\label{th1}
If $\ {\Gamma}/{\beta} + {r^*}/{2} > 1 \ $ the Markov process $Z_0(t)$ is positive recurrent
(ergodic).
\end{theorem}
{\it Proof.} The proof is based on the construction of the Lyapunov function $f$ on the set of all state
of Markov chain, such that the process $f(Z_0(t))$ will be super-martingale. To this aim we use
the criterium of ergodicity from \cite{MP} (Theorem 1.7): In term of our process we have to find a
non-negative function (\textit{Lyapunov function}) $f(x,y)$ such that
$$
(\mathbf{\mathcal{L}_{m-f}} \ f)\left( x,y \right) \le -\varepsilon,
$$
for some $\varepsilon>0$ and for all $(x,y)$ which do not belong to some finite subset $A$
(one must provide it) of the set of all state of the chain.

\begin{figure}[ht!]
\begin{center}
\includegraphics[width=0.8\linewidth]{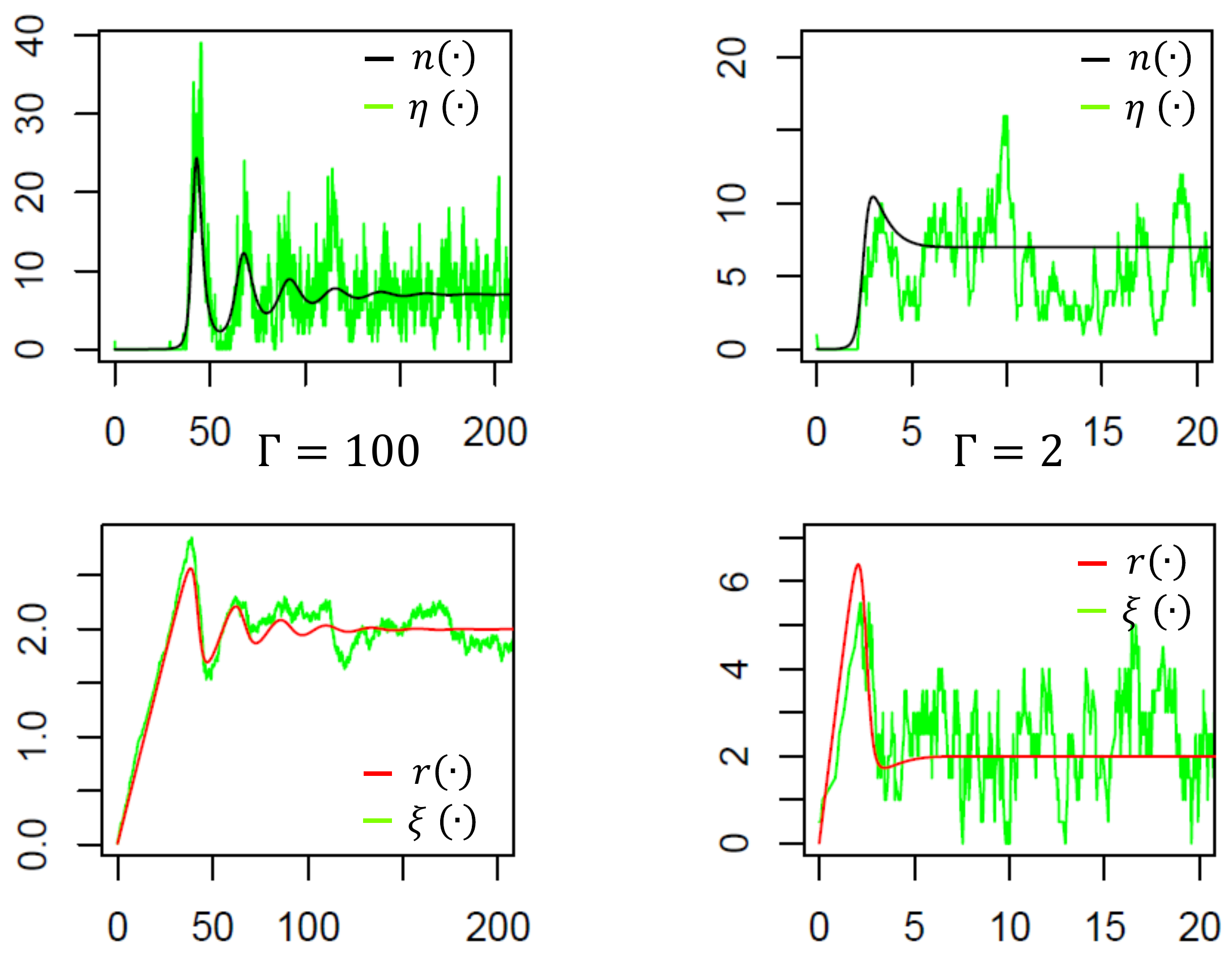}
\end{center}
\caption {Dynamical system and mean-field trajectory: on the top, solutions for $n$ (black) $\Gamma=100$
and for $n$ (black) $\Gamma=2$; on the bottom, solutions for $r$ (red) $\Gamma=100$ and for $r$ (red) $\Gamma=2$.
All starting at $0,01$. Single Markovian trajectories (green) in the \textit{mean-field}
approximation. On the top trajectory samples $\eta$ (green) for $\Gamma=100$ (left) and
$\Gamma=2$ (right); on the bottom trajectory samples $\xi$ (green) for $\Gamma=100$ (left) and $\Gamma=2$ (right).
All for parameters $\alpha=0,01, P=7, \beta=1$.
}%
\label{fig3}
\end{figure}

To this end let us define for our chain the following Lyapunov function :
\begin{equation}\label{Lf1}
f(x,y) = (\Gamma+1) x + y, \ \ \ x\in \frac{1}{\Gamma}\mathbb{Z}_+ , \ y\in \mathbb{Z}_+.
\end{equation}
Then, applying the generator $\mathcal L_{m-f}$ to the Lyapunov function we obtain
\begin{eqnarray*}
(\mathbf{\mathcal{L}_{m-f}} \ f)\left( x,y \right) &= & - \frac{1}{\Gamma} \Bigl( \frac{x n^* +
y r^*}{2} + \alpha x \Bigr) +  \frac{1}{\Gamma} y + \frac{\Gamma+1}{\Gamma} p - \frac{y}{\beta}
\label{eq:3.3gen-0} \\
& =& - \frac{1}{\Gamma} \Bigl( x \Bigl(\frac{n^*}{2}  +  \alpha  \Bigr) + y \Bigl( \frac{\Gamma}{\beta}
+ \frac{r^*}{2} - 1 \Bigr) \Bigr) + \frac{\Gamma+1}{\Gamma} p  \le - \varepsilon \ ,
\label{Lfc1}
\end{eqnarray*}
for some $\varepsilon >0$, when $xn^* + yr^*$ is sufficiently large. This provides the finite
subset $A$, which is defined by
$$
A:= \left\{ (x,y)\in \frac{1}{\Gamma}\mathbb{Z}_+ \times \mathbb{Z}_+ :  \frac{xn^* +
yr^*}{2} + \alpha x + \Bigl( \frac{\Gamma}{\beta} - 1 \Bigr)y
\le (\Gamma +1)p + \Gamma\varepsilon \right\},
$$
and consequently, completes the proof. \hfill  $\Box$
\begin{remark}\label{rem:3.3.0} We note that ${\Gamma}/{\beta} + {r^*}/{2} > 1$ is not a necessary
condition. This choice of parameters makes the construction of Lyapunov function simple. We guess
that the theorem holds true for any choice of positive parameters, $\Gamma, \beta$ and that it may
be proved for other Lyapunov functions.
\end{remark}
Behaviour of a single trajectories in the mean-field approximation is illustrated by Figure \ref{fig3}.
Note that qualitatively it is more \textit{wiggling} then the arithmetic mean for $N=50$, but less
\textit{spiking} then that for $N=10$, see Figures \ref{fig1} and \ref{fig2}.

Taking into account the relation between representation (\ref{s1}) for DS and the form
of generator (\ref{eq:3.3gen-0}) we deduce that corresponding to mean-field approximation
dynamical system gets the form
\begin{equation}\label{eq:3.3s1-0}
\partial_t \binom{{r}}{{n}} = \Bigl( \frac{r n^* + n r^*}{2} +
\alpha r \Bigr) \binom{-1/\Gamma}{1} + n \binom{1/\Gamma}{-1} + p \binom{1/\Gamma}{0} +
\beta^{-1} n \binom{0}{-1}  \ .
\end{equation}
{Then, in the stationary regime, we expect that its stationary points coincides with stationary points of DS.
Indeed, the following theorem holds.}

\begin{theorem}\label{st.points1}
In the stationary regime, when $\partial_t {r} = \partial_t {n} = 0$,
(\ref{eq:3.3s1-0}) defines the stationary expectation values $(r^{**}, n^{**})$ of the
process $Z_0(t)=(\xi_0(t),\eta_0(t))$ which coincides with stationary points of DS $(r^{*}, n^{*})$:
$(r^{**}, n^{**}) = (r^{*}, n^{*}).$
\end{theorem}

 {\it Proof.} Since by (\ref{s2})
\begin{equation*}
\left\{ \begin{array}{rcl}
r^* & = & p(1 + \beta)/(\beta p + \alpha) \ , \\
n^*  &=& \beta p \ ,
\end{array} \right.
\end{equation*}
the explicit calculations with help of (\ref{eq:3.3s1-0}) yield $r^{**} = r^*$ and $n^{**} = n^*$. $\Box$

This confirms (at least for stationary regime) a consistency of expectations of
$Z_0(t)=(\xi_0(t),\eta_0(t))$ for the process (\ref{eq:3.3gen-0}) with the values of $r^*$ and
$n^*$ taken for the mean-field generator from DS (\ref{DS1}), (\ref{DS2}).

\subsection{One-unit Markov process and spikes}\label{sec.3.4}

In a sense we would like to study the case, which is intermediate between global ($N \gg 1$) and
one-unit ($N=1$) limits.
The one-unit Markov process describes the system for $N=1$. The corresponding generator
(\ref{Mark-gen}) is simply $\mathbf{L}_{N=1}$:
\begin{eqnarray}\label{Mark-gen-one}
(\mathbf{L_1} g)\left( x,y \right) &=& x \, y \ [g(x -\frac{1}{\Gamma}, y + 1)- g(x,y)] \\
&+& \alpha x \ [g(x -\frac{1}{\Gamma}, y + 1)- g(x,y)] \nonumber \\
&+& y \ [g(x+\frac{1}{\Gamma}, y - 1)- g(x,y)] \nonumber \\
&+& p \ [g(x +\frac{1}{\Gamma}, y)- g(x,y)] \nonumber \\
&+& \beta^{-1} y \ [g(x, y - 1)- g(x,y)] \ , \nonumber
\end{eqnarray}
where we put $(r_1 , n_1) = (x , y) \in \frac{1}{\Gamma}\mathbb{Z}_+ \times \mathbb{Z}_+ $
\begin{figure}[ht!]
\begin{center}
\includegraphics[width=0.8\linewidth]{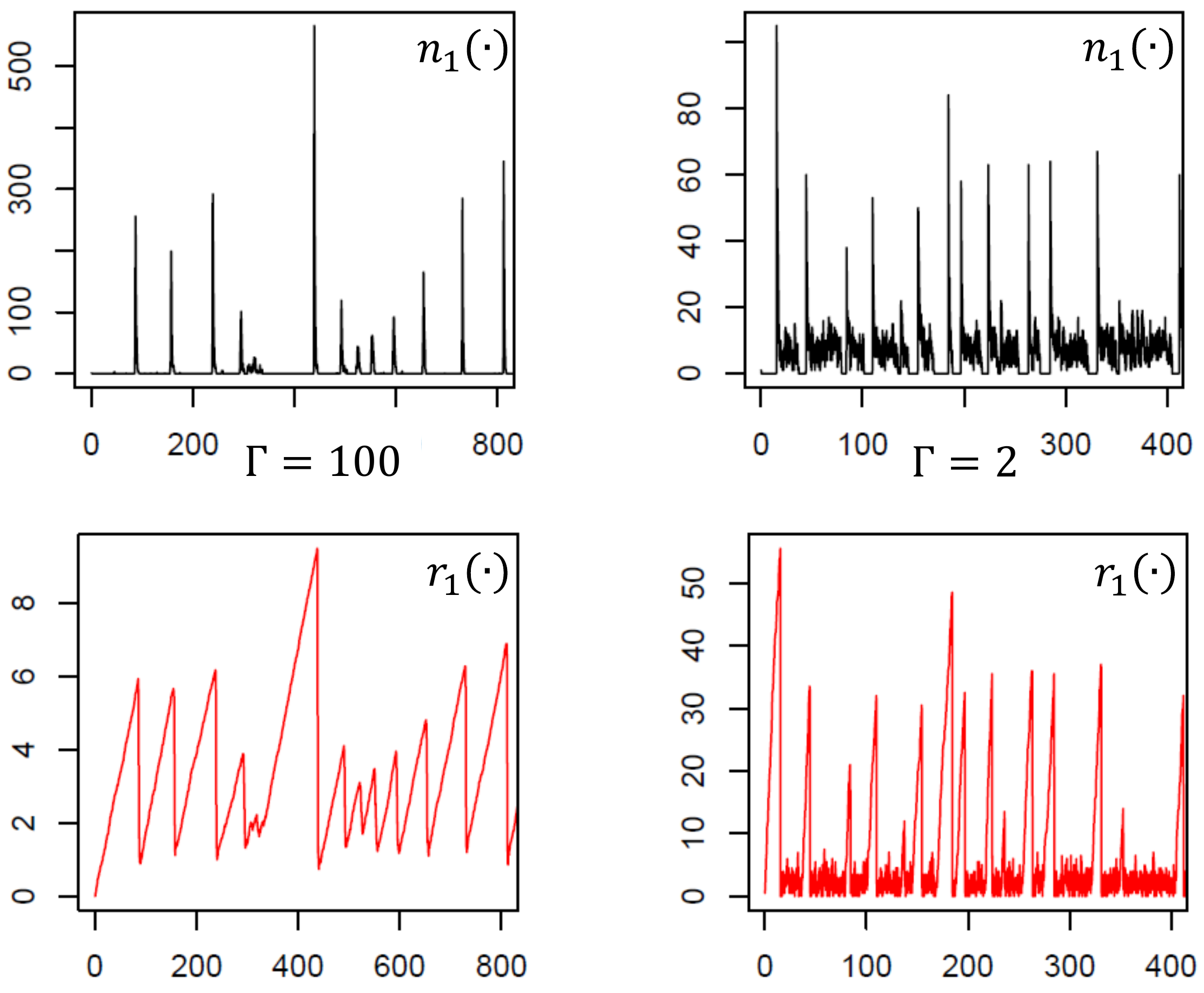}
\end{center}
\caption {Markov processes for N=1. Top-left, bottom-left,
the trajectory samples of stochastic Markov processes:  {$n_1$} (black) and {$r_1$} (red),
for the parameters $\alpha=0,01, \Gamma=100, P=7$ and $\beta=1$.
Similarly, {$n_1$} (black) and {$r_1$} (red) on the top-right,
bottom-right corresponds to parameters:  $\alpha=0,01, \Gamma=2, P=7$ and $\beta=1$.}
\label{fig5}
\end{figure}
In this case the normalising parameter $N=1$ and we have only \textit{single} Markov trajectory
(\ref{eq:3.12}) $X_{N=1}(t){=\binom{r_1(t)}{n_1(t)}}$. Therefore, this is \textit{not} a {density dependent}
Markov process \cite{EK} driven by a dynamical system for densities.

According to (\ref{Mark-gen-one}) the first component {$r_1(t)$} $\in ({1}/{\Gamma})\mathbb{Z}_+$
of $X_{N=1}(t)$  is the random inversion coefficient, which is jumping
between \textit{zero} and \textit{infinity} with increment $1/\Gamma$. The higher is excited level,
the larger is the coefficient of inversion $\xi$ in the system. It is zero for the ground state of
the system. On the other hand, the second component {$n_1(t)$} $\in \mathbb{Z}_+$ of $X_{N=1}(t)$
is counting the instant number of photons in the system. It is also jumping
between \textit{zero} and \textit{infinity} but with the increment \textit{one}.

We note that interpretation of our model for $N=1$ is straightforward and perfectly understandable
in the framework of Remark \ref{rem:3.1}, when the total number: non-excited $+$ exited atoms
is \textit{not fixed}. But in fact it is also equivalent to a model with a \textit{single} atom $N=1$.
Indeed, since the coefficients of transition intensities in $\mathbf{L}_1$ are motivated by Einstein
equations of Section \ref{sec.2.1}, the \textit{one-unit} system for $N=1$ allows interpretation as the
model of a \textit{single} \textit{infinitely-many} level atom with the spacing between levels
equals to $1/\Gamma$. This atom has the instant inversion coefficient {$r_1$} and it is embedded
into the random photon environment with the instant photon intensity {$n_1$}.

We illustrate the behaviour of the one-unit system on Figure \ref{fig5}.
\begin{theorem}\label{t2}
When $\Gamma \ge \beta$ the Markov chain, which is governed by generator $\mathbf{L_1}$ is
positive recurrent Markov chain.
\end{theorem}
{\it Proof.} The proof follows the same line of reasoning as above for individual component
in mean-field. Moreover, the same Lyapunov function (\ref{Lf1}) can be applied in this case. Indeed,
using  (\ref{Lfc1}) again we obtain
\begin{equation}\label{Lfc2}
(\mathbf{L_1} f)\left( x,y \right) =  - \frac{1}{\Gamma} \Bigl( xy + \alpha x +
\Bigl( \frac{\Gamma}{\beta} - 1 \Bigr)y \Bigr) + \frac{\Gamma+1}{\Gamma} p  \le - \varepsilon
\end{equation}
for some $\varepsilon >0$ and $x+y$ large enough, and when $\Gamma > \beta$. This also provides the set
$A$ defined by:
$$
A = \left\{ (x,y): \ xy + \alpha x + \Bigl( \frac{\Gamma}{\beta} - 1 \Bigr)y \le (\Gamma +1)p +
\Gamma \varepsilon \right\}.
$$
This gives the proof of assertion. \hfill  $\Box$
\subsection{One-unit process: statistics of spikes}\label{sec.3.5}

{The main characteristics of trajectory for $n_1(\cdot)$ in one-unit model is the presence of well featured
spikes and plateaus just before spikes. Then a natural question concerns distributions of the spikes amplitude
and of the lengths of plateau as well as their correlation. Since to obtain these distributions and correlation
explicitly (analytically) is a difficult problem, we propose only some results of numerical statistical
analysis.}
\begin{figure}[ht!]
\begin{center}
\includegraphics[width=0.9\linewidth]{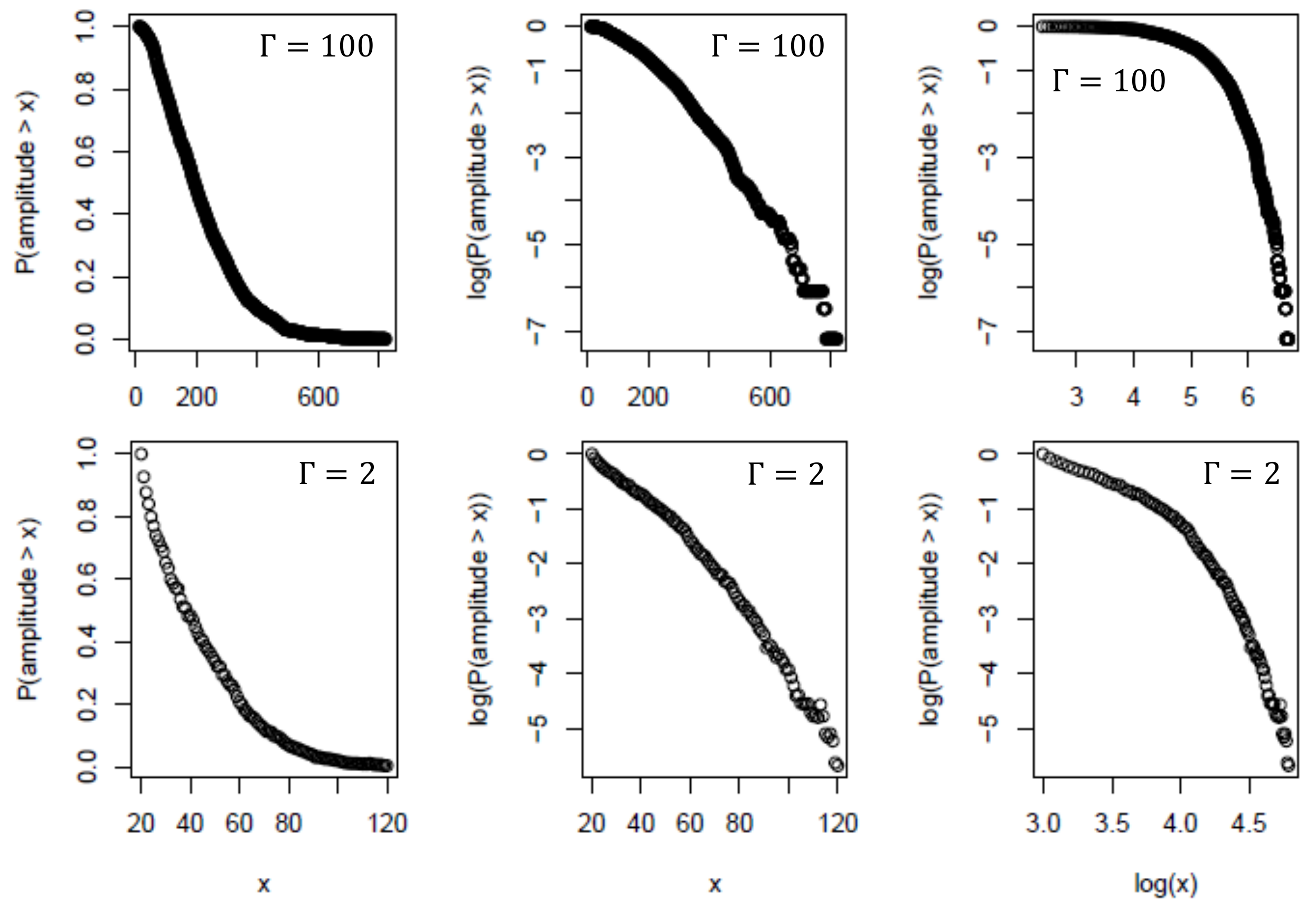}
\end{center}
\caption {Spikes amplitude statistics: Markov process for N=1, $\Gamma=100$ and $\Gamma=2$.
Simulations were performed with the parameters:  $\alpha=0,01, P=7$ and $\beta=1$. We see that the estimated
tail probability follows the exponential law (middle column). }
\label{fig6}
\end{figure}
\begin{figure}[ht!]
\begin{center}
\includegraphics[width=0.8\linewidth]{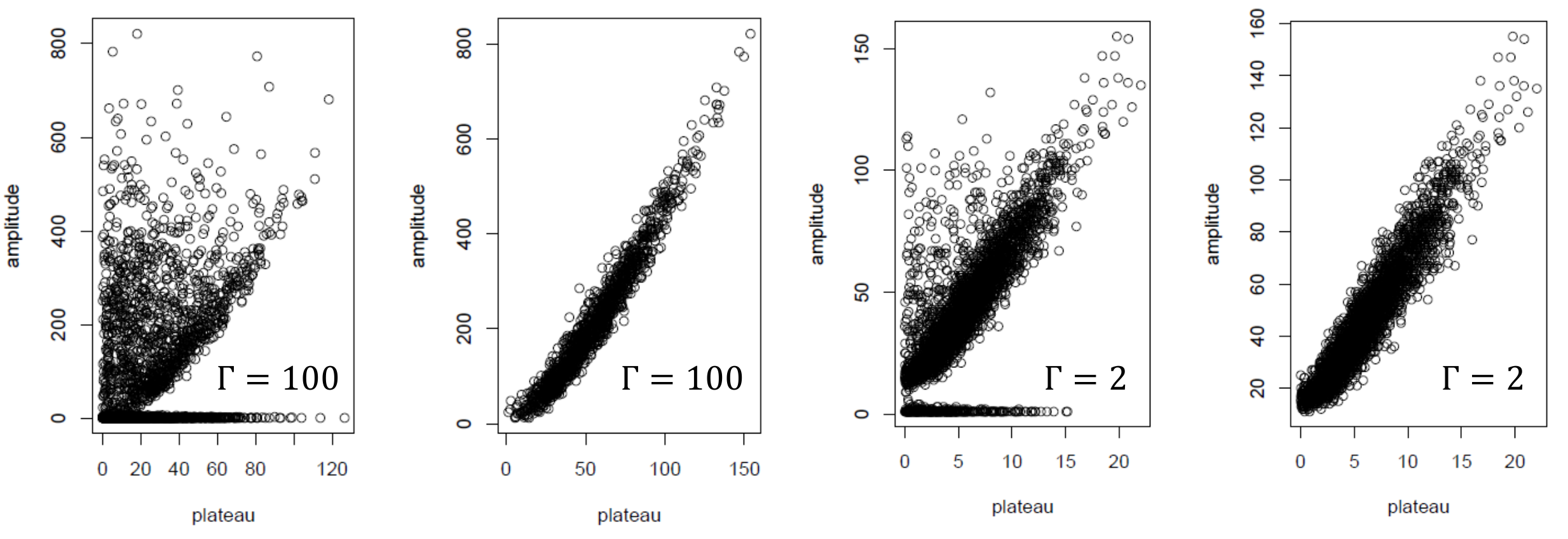}
\end{center}
\caption {Plateau versus amplitude statistics: Markov process for one-unit process, $N=1$,
$\Gamma=100$ and $\Gamma=2$.
On the first and third scatterplots (for both $\Gamma=100$ and $\Gamma=2$) the plateau is defined
as interval when $n_1(t)=0$. On the second and forth scatterplots the plateau is defined as the
time interval when $n_1(t)\le 10$ for the respectively ($\Gamma=100$ and $\Gamma=2$) same simulated
trajectory. Other parameters in simulations were, as before, $\alpha=0,01, P=7$ and $\beta=1$.
 }%
\label{fig7}
\end{figure}
The statistics of spikes' \textit{amplitudes} for the one-unit Markov process with generator $\mathbf{L}_1$
is presented on the next Figure \ref{fig6}. On the simulated trajectory for a fixed value $a$ we calculate
the frequency of times when the spikes are greater than $a$. For the case $\Gamma=100$ the minimal threshold
for amplitude of spikes was chosen $a_0=10$, thus, we estimate the probability of the amplitude of spike greater
than $a$ \textit{given} the amplitude is greater then $10$. For $\Gamma=2$ the minimal threshold is $20$. These
choice of minimal thresholds makes the graphs more \textit{readable}. We observe in Figure~\ref{fig6} that
the tail distribution of the amplitude in both cases ($\Gamma=100$ and $\Gamma=2$) is similar to exponential
(median graphs).

According the typical trajectories for the case when $\Gamma$ is large, see Figure~\ref{fig5} for $\Gamma=100$
we observe the presence of plateau (a time interval when $\eta(t)=0$). We expect that
the  plateau interval and a successive amplitude will be \textit{positively correlated}. The first
scatterplot, see the left hand-side scatterplot on Figure~\ref{fig7} gives us an idea about
\textit{positive correlation} between plateau and amplitude. But the figure also indicate the
presence of large plateau with very small successive amplitude. It provide the new definition of plateau: a
plateau is the time interval when $\eta(t) \le thr$, where parameter $thr$ is threshold parameter. With new
definition of plateau the tendency of positive correlation is more obvious, and moreover the dependency is not
linear, which is illustrated for the case $\Gamma=100$ and $thr=10$ in right scatterplot on Figure~\ref{fig7}.

It is instructive to compare the \textit{one-unit} Markov process with generator $\mathbf{L}_1$
(\ref{Mark-gen-one}) with the \textit{single} Markov trajectory in the mean-field approximation with
generator $\mathbf{\mathcal{L}}$ (\ref{Mark-gen-single}). Note that only intensities of transitions in the
first term are different. For the mean-field approximation the intensity is more smooth than for the
one-unit Markov process, cf. the corresponding pictures Figure~\ref{fig3} and Figure~\ref{fig5}.

\vspace{0.5cm}

\textbf{Acknowledgements}

This work was supported by Funda\c{c}\~ao de Amparo \`a Pesquisa do Estado de S\~ao Paulo (FAPESP) under Grant 2016/25077-4.

AY also thanks Conselho Nacional de Desenvolvimento Cient\'ifico e Tecnol\^ogico (CNPq) grant 301050/2016-3 and
FAPESP grant 2017/10555-0.

VAZ is grateful to Instituto de Matem\'{a}tica e Estat\'{\i}stica of the University of S\~{a}o Paulo and to
Anatoly Yambartsev for a warm hospitality. His visits were supported FAPESP grant 2016/25077-4.



\end{document}